\documentclass[10pt]{article}
\usepackage{amsmath}
\usepackage{latexsym}
\usepackage{amssymb}
\usepackage{amsfonts}
\usepackage{amsthm}
\title{HOW MANY STRUCTURE CONSTANTS \\
DO EXIST IN RIEMANNIAN GEOMETRY}
\author{  J.-F. POMMARET  \\
ORCID:0000-0003-0907-2601}
\date{  }
\begin{document}
\maketitle

\noindent
{\bf ABSTRACT}  \\

After reading such a question, any mathematician will say that, according to a well known result of L.P. Eisenhart found in 1926, the answer is " One " of course, namely the only constant allowing to describe the so-called " {\it constant riemannian curvature} " condition. The purpose of this paper is to prove the contrary by studying the case of two dimensional riemannian geometry in the light of an old work of E. Vessiot published in 1903 but {\it still totally unknown today} after more than a century. In fact, we shall compute locally the {\it Vessiot structure equations} and prove that there are indeed " Two " {\it Vessiot structure constants} satisfying a single {\it linear Jacobi condition} showing that one of them must vanish while the other one must be equal to the known one. This result depends on deep mathematical reasons in the formal theory of Lie pseudogroups, which are involving both the Spencer $\delta$-cohomology and diagram chasing in homological algebra. Another similar example will illustrate and justify this comment out of the classical tensorial framework of the famous " {\it equivalence problem} ". The case of contact trasformations will also be studied. Though it is quite unexpected, we shall reach the conclusion that the mathematical foundations of both classical and conformal riemannian geometry must be revisited. We have treated the case of conformal geometry in a recent arXiv preprint. \\

\vspace{4cm}

\noindent
{\bf KEY WORDS}:      \\

\noindent
Lie pseudogroups; Formal integrability; Vessiot structure equations; Riemannian geometry; \\
Equivalence problem.   \\

\newpage

\noindent
{\bf 1) INTRODUCTION}  \\
When $X$ is a manifold of dimension $n$ with local coordinates $(x^1, ..., x^n)$, we first sketch the discovery of Vessiot ([8],[21]) {\it still not known today after more than a century} for reasons which are not scientific at all ([10]). Roughly, using standard notations of jet theory, a Lie pseudogroup $\Gamma \subset aut(X)$ is made by finite invertible transformations $y=f(x)$ solutions of a (nonlinear in general) system ${\cal{R}}_q\subset {\Pi}_q$ while, using vertical bundles, the infinitesimal transformations $\xi\in \Theta$ are solutions of the linearized system $R_q=id_q^{-1}(V({\cal{R}}_q))\subset J_q(T)$ where $T=id^{-1}(V(X\times X)$ is the tangent bundle of $X$, $id:X \rightarrow X \times X:x \rightarrow (x,x)$ the identity map and $id_q=j_q(id)$ the $q$-jet of the identity. When $\Gamma$ is transitive, there is a canonical epimorphism ${\pi}^q_0:R_q\rightarrow T$. Also, as changes of source $x$ commute with changes of target $y$, they exchange between themselves any generating set of {\it differential invariants} $\{{\Phi}^{\tau}(y_q)\}$ of order $q$. Then one can introduce a {\it natural bundle} $\cal{F}$ over $X$, also called {\it bundle of geomeric objects}, by patching changes of coordinates of the form $\bar{x}=f(x), \bar{u}=\lambda(u,j_q(f(x))$ thus obtained (See examples below). A section $\omega$ of $\cal{F}$ is called a {\it geometric object} or {\it structure} on $X$ and transforms like ${\bar{\omega}}(f(x))=\lambda(\omega(x),j_q(f)(x))$ or simply $\bar{\omega}=j_q(f)^{-1}(\omega)$. This is a way to generalize vectors and tensors ($q=1$) or even connections ($q=2$). As a byproduct we have $\Gamma=\{f\in aut(X){\mid} j_q(f)^{-1}(\omega)=\omega\}$ and we may say that $\Gamma$ {\it preserves} $\omega$. Replacing $j_q(f)$ by $f_q$, we also obtain the {\it Lie form} ${\cal{R}}_q=\{f_q\in {\Pi}_q{\mid} f_q^{-1}(\omega)=\omega\}$. Coming back to the infinitesimal point of view and setting $f_t=exp(t\xi)\in aut(X), \forall \xi\in T$, we may define the {\it ordinary Lie derivative} with value in the vector bundle $F_0={\omega}^{-1}(V({\cal{F}}))$ by the formula :\\
\[ {\cal{D}}\xi={\cal{L}}(\xi)\omega=\frac{d}{dt}j_q(f_t)^{-1}(\omega){\mid}_{t=0} \Rightarrow \Theta=\{\xi\in T{\mid}{\cal{L}}(\xi)\omega=0\}  \]
and we say that $\cal{D}$ is a {\it Lie operator} because ${\cal{D}}\xi=0,{\cal{D}}\eta=0\Rightarrow {\cal{D}}[\xi,\eta]=0$ as we already saw.\\
In the jet framework at any order $q$, we shall introduce by linearity as in ([8],[11],[14],[15]) the {\it formal Lie derivative} $L({\xi}_q)$ in such a way that ${\cal{L}}(\xi)=L(j_q(\xi))$. It follows thet the infinitesimal Lie equations defining $R_q$ can be written in the so-called {\it Medolaghi form} $L({\xi}_q)\omega=0$ with coefficients depending on $j_1(\omega)$ in a very specific fashion ([6],[8]):   \\
\[   {\Omega}^{\tau}\equiv (L({\xi}_q)\omega)^{\tau}= - L^{\tau \mu}_k(\omega (x)){\xi}^k_{\mu} + {\xi}^r {\partial}_r {\omega}^{\tau}(x) = 0  \]
Let us suppose that the symbol $g_q\subset S_qT^* \otimes T$ is involutive, in such a way that this system becomes formally integrable and thus involutive, that is all the equations of order $q+r$ could be obtained by differentiating $r$ times {\it only}, $\forall r\geq 0$. Then, as we shall see in the following examples, $\omega$ must satisfy certain (non-linear in general) integrability conditions of the form: \\
\[    I(j_1(\omega)) =c (\omega)     \]
called {\it Vessiot structure equations}, depending on a certain number of {\it Vessiot structure constants} $c$ eventually satisfying algebraic {\it Jacobi conditions} $J(c)=0$ and we let the reader compare this situation to the Riemannian or contact cases ([14]). In the second section, we shall treat a specific example when $n=2$ while in the third section, we shall treat the Riemannian case with full details when $n=2$. The most striking result of this paper is that, though at first sight there does not seem to be any link between these two examples, we shall discover at the end of the paper that they are in fact ... identical !. We want to point out that these structure equations were perfectly known by E. Cartan (1869-1951) who {\it never} said that these results were at least competing with or even superseding the corresponding {\it Cartan structure equations } that he has developed about at the same time for similar purposes ([1]). The underlying reason is of a purely personal origin related to the {\it differential Galois Theory} within a kind of "{\it mathematical affair} " involving the best french mathematicians of that time ([9]). The original letters, given to the author of this paper by M. Janet, a friend of E. Vessiot, have ben published in ([10]) and have been put as a deposit in the main library of Ecole Normale Sup\'{e}rieure in Paris for future historical studies. \\
             Finally, we can choose for the generating CC ${\cal{D}}_1$ of ${\cal{D}}$ the first order linearization of a non-linear version described by the Vessiot structure equations:  \\
\[     \frac{\partial I}{\partial j_1(\omega)}(j_1(\omega)) j_1(\Omega) =\frac{\partial c}{\partial \omega}(\omega)\Omega   \]
that is {\it exactly} what is usually done for the flat Minkowski metric in general relativity ([12],[16],[17]).  \\

\newpage

\noindent
{\bf 2) MOTIVATING EXAMPLES}   \\
We show that the Vessiot structure equations may even exist when $n=1$ ([12]). In the remaining of this paper, the reader may refer to ([7],[19] or [11]) for the elements of homological algebra allowing to chase in the commutative diagrams that we shall present.   \\

\noindent
{\bf EXAMPLE 2.1}:  While the {\it affine} transformations $y=ax+b$ are solutions of the second order linear system $y_{xx}=0$, the sections of the corresponding linearized systems are respectively satisfying ${\xi}_{xx}=0$. The only generating differential invariant $\Phi\equiv y_{xx}/y_x$ of the affine case transforms like $u=\bar{u}{\partial}_xf+({\partial}_{xx}f/{\partial}_xf)$ when $\bar{x}=f(x)$. The corresponding geometric object defined by the section $u=\gamma (x)$ does transform like the Christoffel symbols, namely:  
\[        \gamma(x)=\bar{\gamma}(f(x)){\partial}_xf + ({\partial}_{xx}f/{\partial}_xf)   \]
For this, if $\gamma$ is the geometric object of the affine group $y=ax+b$ and $0\neq \alpha=\alpha (x)dx \in T^*$ is a $1$-form, we consider the object $\omega=(\alpha,\gamma)$ and get at once the two {\it Medolaghi equations}:\\  
\[  {\cal{L}}(\xi)\alpha\equiv \alpha {\partial}_x\xi + \xi {\partial}_x\alpha =0, \hspace{1cm} {\cal{L}}(\xi)\gamma\equiv {\partial}_{xx}\xi+\gamma {\partial}_x\xi+ \xi {\partial}_x\gamma =0  \]
Differentiating the first equation and substituting the second, we get the zero order equation:  \\
\[  \xi (\alpha {\partial}_{xx}\alpha-2({\partial}_x\alpha)^2+\alpha \gamma {\partial}_x\alpha-{\alpha}^2{\partial}_x\gamma)=0\hspace{5mm} \Leftrightarrow \hspace{5mm} \xi {\partial}_x(\frac{{\partial}_x\alpha}{{\alpha}^2} - \frac{\gamma}{\alpha} )=0  \]
and the {\it Vessiot structure equation} ${\partial}_x\alpha-\gamma \alpha=c{\alpha}^2$. Alternatively, setting $\beta=-1/\alpha\in T$, we get ${\partial}_x\beta+\gamma\beta=c$. With $\alpha=1, \beta=-1, \gamma=0 \Rightarrow c=0$ we get the translation subgroup $y=x+b$ while, with $\alpha=1/x, \beta=-x, \gamma=0 \Rightarrow c=-1$ we get the dilatation subgroup $y=ax$. \\
Working now with isometries, we just need to set $\omega = {\alpha}^2$ in order to obtain the Killing equation $ 2 \omega
 {\xi}_x + \xi {\partial}_x \omega =0$ and the corresponding Vessiot structure equation 
 ${\partial}_x \omega - 2 \omega \gamma = c' {\omega}^{\frac{3}{2}} $ with $c'=2c$.  \\
Similarly, if $\nu$ is the geometric object of the projective group $y=(ax+b)/(cx+d)$, transforming like the well known {\it Schwarzian differential invariant} $\Psi=(y_{xxx}/y_x) - \frac{3}{2} (y_{xx}/y_x)^2=d_x\Phi - \frac{1}{2} {\Phi}^2$, we may consider the new geometric object $\omega=(\gamma,\nu)$ and get the only {\it Vessiot structure equation} ${\partial}_x\gamma - \frac{1}{2}{\gamma}^2-\nu=0$ in a coherent way, without any structure constant ([18]).\\

\noindent
{\bf EXAMPLE 2.2}:  With $m=n=2$, let us now consider the Lie group of transformations $\{ y^1=ax^1 +b, y^2=cx^2 +d \mid a,b,c,d=cst, ac=1 \}$ as an algebraic Lie pseudogroup $\Gamma$. It is easy to exhibit the corresponding first order system ${\cal{R}}_1$ of finite Lie equations in Lie form by introducing the three generating differential invariants and the corresponding Lie form:  \\
\[  {\Phi}^1 \equiv \frac{y^1_2}{y^1_1}=0, \,\,\,   {\Phi}^2\equiv \frac{y^2_1}{y^2_2}=0, \,\,\,  {\Phi}^3 \equiv  y^1_1y^2_2=1   \]
The details of the corresponding tricky computations, first done in 1978 ([8]), have been improved in 2016 ([14]) and are again revisited in this paper. As we shall see, its major interest is to work out the two Vessiot structure constant existing like for the Riemannian structure but without having {\it any} tensorial framework ([5]).  \\
First of all, we notice that this system is finite type with a vanishing second order symbol and is thus formally integrable but not involutive. For this, we may introduce the six generating differential invariants obtained after one prolongation, exactly like the six Christoffel symbols in Riemannian geometry:  \\
\[   {\Phi}^4 \equiv \frac{y^1_{11}}{y^1_1}=0, {\Phi}^5 \equiv \frac{y^1_{12}}{y^1_1}=0, {\Phi}^6\equiv \frac{y^1_{22}}{y^1_1}=0, {\Phi}^7 \equiv \frac{y^2_{22}}{y^2_2} = 0, {\Phi}^8 \equiv \frac{y^2_{12}}{y^2_2} =0, {\Phi}^9 \equiv \frac{y^2_{11}}{y^2_2}=0   \]
while noticing that ${\Phi}^3(1-{\Phi}^1{\Phi}^2) \equiv y^1_1 y^2_2 - y^1_2 y^2_1\neq 0$ is changing like the Jacobian. \\
Looking to the way these invariants are transformed under an arbitrary change $ \bar{x}=\varphi (x)$, we obtain for example $u^1 = ({\partial}_2{\varphi}^1 + {\bar{u}}^1{\partial}_2{\varphi}^2) / ({\partial}_1{\varphi}^1 + {\bar{u}}^1{\partial}_1 {\varphi}^2)$ and so on for describing the natural fiber bundle ${\cal{F}}$ with local coordinates $(x^1,x^2; u^1,u^2,u^3)$ and section $\omega = ({\omega}^1,{\omega}^2, {\omega}^3)$ becoming $(0,0,1)$ with our choice and ${\omega}^3(1- {\omega}^1{\omega}^2)=1 \neq 0$. Passing to the infinitesimal point of view, we obtain the first order system $R_1 \subset J_1(T)$ in the Medolaghi form with jet notation $\Omega \equiv L({\xi}_1) \omega=0$ as in ([8]) and may compare it to the three equations of the Killing system in dimension $n=2$:  \\
\[   \left\{  \begin{array}{rcl}
{\Omega}^1 & \equiv & {\xi}^1_2 + {\omega}^1 {\xi}^2_2 - {\omega}^1 {\xi}^1_1 - ({\omega}^1)^2 {\xi}^2_1 + {\xi}^r {\partial}_r {\omega}^1= 0  \\
{\Omega}^2 & \equiv & {\xi}^2_1 + \omega^2 {\xi}^1_1 - {\omega}^2{\xi}^2_2 - ({\omega}^2)^2{\xi}^1_2 
+ {\xi}^r {\partial}_r {\omega}^2 =0  \\
 {\Omega}^3 & \equiv & {\omega}^3 ({\xi}^1_1 + {\xi}^2_2) + {\omega}^1{\omega}^3 {\xi}^2_1 +{\omega}^2{\omega}^3 {\xi}^1_2  + {\xi}^r {\partial}_r {\omega}^3 =0
 \end{array}  \right.   \]
 that we can extend to six {\it intermediate} new equations like in ([14]), including in particular:  \\
 \[  {\Omega}^4 \equiv  {\xi}^1_{11}+ {\omega}^1 {\xi}^2_{11} + {\omega}^4{\xi}^1_1  + 2({\omega}^5 - {\omega}^1 {\omega}^4){\xi}^2_1 + {\xi}^r {\partial}_r {\omega}^4 =0 \]
 \[ {\Omega}^5 \equiv {\xi}^1_{12} + {\omega}^1{\xi}^2_{12} + {\omega}^4 {\xi}^1_2 = ({\omega}^6 - {\omega}^1{\omega}^5){\xi}^2_1 + {\omega}^5{\xi}^2_2 + {\xi}^r {\partial}_r{\omega}^5=0  \]
  \[   {\Omega}^6 \equiv {\xi}^1_{22} + {\omega}^1 {\xi}^2_{22}++ 2 {\omega}^6 {\xi}^2_2 +  2 {\omega}^5 {\xi}^1_2 - {\omega}^6 {\xi}^1_1 - {\omega}^1 {\omega}^6 {\xi}^2_1 + {\xi}^r {\partial}_r{\omega}^6 = 0  \]
Taking into account these new invariants bringing for example six relations like:   \\
\[  {\partial}_1{\omega}^1 - {\omega}^5 + {\omega}^1 {\omega}^4 = 0, \hspace{1cm} {\partial}_2 {\omega}^1 - {\omega}^6 + {\omega}^1 {\omega}^5 = 0  \]
\[   {\partial}_1 {\omega}^2 - {\omega}^9  + {\omega}^2 {\omega}^8 = 0 , \hspace{1cm}  {\partial}_2 {\omega}^2 - {\omega}^8 + {\omega}^2 {\omega}^7 =0   \]
\[  {\partial}_1 {\omega}^3 - {\omega}^3 ( {\omega}^4 + {\omega}^8)=0, \hspace{1cm} {\partial}_2 {\omega}^3 - {\omega}^3 ({\omega}^5 + {\omega}^7) = 0  \]
The determinant of the $6 \times 6$ matrix with respect to $({\omega}^4, ... , {\omega}^9)$ is ${\omega}^3 ( 1 - {\omega}^1{\omega}^2) \neq 0$. \\
After tedious but elementary substitutions, we obtain for example:  \\
\[   d_2{\Omega}^4 - d_1 {\Omega}^5 \equiv ({\partial}_2 {\omega}^4 - {\partial}_1{\omega}^5)({\xi}^1_1 +{\xi}^2_2) + {\xi}^r{\partial}_r ({\partial}_2 {\omega}^4 - {\partial}_1{\omega}^5) = 0  \]
However, we have also:  \\
\[  {\omega}^3 (1 - {\omega}^1{\omega}^2) ({\xi}^1_1 + {\xi}^2_2) + {\xi}^r {\partial}_r ({\omega}^3 (1 - {\omega}^1{\omega}^2))=0  \]
Replacing $({\omega}^4, ... ,{\omega}^9)$ by their rational expressions in $j_1(\omega)$, the quotient of ${\partial}_2{\omega}^4 - {\partial}_1 {\omega}^5$ by ${\omega}^3(1 - {\omega}^1{\omega}^2)$ is well defined, say equal to $c(x)$, and we obtain the new zero order equation ${\xi}^r{\partial}_r c(x)=0$ contradicting the formal integrability of the given system unless we obtain the two Vessiot structure equations with the Vessiot structure constant $c', c''$, namely:  \\
\[         {\partial}_2 {\omega}^4 - {\partial}_1 {\omega}^5 = c' \,\, {\omega}^3(1- {\omega}^1{\omega}^2), \hspace{1cm}
{\partial}_1 {\omega}^7 - {\partial}_2 {\omega}^8 = c'' \,\, {\omega}^3 (1 - {\omega}^1 {\omega}^2)   \]
It remains to prove that there is only one Vessiot structure equation of order two with a single structure constant. For this, first of all we notice that:  \\
\[   {\partial}_2({\omega}^4 + {\omega}^8) - {\partial}_1 ({\omega}^5 + {\omega}^7) = 0 \,\,  \Rightarrow \,\,   c'-c'' = 0 \,\,  \Rightarrow \,\,  c'=c''=c   \]
Now, if $R_q \subset J_q(T)$, let us define $J^0_q(T)$ by the short exact sequence:  \\
\[     0 \rightarrow J^0_q(T) \rightarrow J_q(T) \stackrel{{\pi}^q_0}{\longrightarrow} T \rightarrow 0  \]
and set $R^0_q=R_q \cap J^0_q(T) \subset J_q(T)$. We have the commutative and exact diagram:  \\
\[      \begin{array}{rcccccl}
      & 0 &  & 0   &  &  &    \\
      & \downarrow &   &  \downarrow  &  &  &  \\
 0 \rightarrow &  g_{q+1} &  =  &  g_{q+1}  &\rightarrow    & 0 &  \\
      &   \downarrow  &  &  \downarrow &  &  \downarrow &    \\
  0 \rightarrow &  R^0_{q+1} &  \rightarrow  & R_{q+1}  &  \rightarrow &  T &       \\
      &  \downarrow &   &  \downarrow &  &  \parallel  &   \\
  0 \rightarrow &  R^0_q  &   \rightarrow &  R_q & \rightarrow & T  &  \rightarrow 0   \\
                         &     &   &   &   &\downarrow    &       \\
                          &   &   &   &   &   0  &   
                        \end{array}                   \]
It follows that $ R_{q+1} \rightarrow R_q$ is an epimorphism if and only if  $R^0_{q+1} \rightarrow R^0_q$ and $R_{q+1} \rightarrow T$ are both epimorphisms. It just remains to use successively $q=1$ and $q=2$.   \\                          
In the present situation, constructing the same diagram as the one used in the study of the Killing system (See [8] or [11] for the details), we have the commutative and exact diagram allowing to construct the second order CC when $g_2=0,g_3=0$, namely:  \\
\[  \begin{array}{rcccccccl}
    &   &    &   0  & &  0  &  &  &     \\
     &  &    &   \downarrow  &  &  \downarrow  &   &  &  \\
   &  0  & \rightarrow &  S_3T^*\otimes T  &  \rightarrow &  S_2T^*\otimes F_1 &  \rightarrow & F_2 & \rightarrow 0   \\
   &   &   &  \downarrow &   &  \downarrow &  &   \\
 &  0  & \rightarrow & T^* \otimes  S_2T^*\otimes T  &  \rightarrow &  T^* \otimes T^*\otimes F_1 &  \rightarrow & 0 &   \\          
&   &   &  \downarrow &   &  \downarrow &  &   \\
0 \rightarrow &  {\wedge}^2T^* \otimes g_1 &  \rightarrow &  {\wedge}^2T^* \otimes T^* \otimes T & \rightarrow &  {\wedge}^2T^* \otimes F_1 & \rightarrow &  0  &   \\
   &   &   &  \downarrow &   &   \downarrow  &  &   &  \\
   &   &  &   0  &  &   0   &   &   &
   \end{array}    \]
We obtain the isomorphism $F_2 \simeq {\wedge}^2T^* \otimes g_1$ by a snake chase and deduce  thus the relation $dim(F_2)= dim ({\wedge}^2T^*\otimes g_1)=dim (g_1)=1$ because $n=2$. Of course, as $dim(F_1)=3$, we also obtain $dim(F_2)= dim(S_2T^*\otimes F_1) - dim(S_3T^*\otimes T)= 9 - 8=1$ in a coherent way. \\
For the sake of completeness, we provide the only component of the second order CC, namely:   \\
\[  {\Omega}^1 \equiv {\xi}^1_2=0, {\Omega}^2 \equiv {\xi}^2_1=0, {\Omega}^3 \equiv {\xi}^1_1 + {\xi}^2_2=0 \,\,\,  \Rightarrow \,\,\, d_{11}{\Omega}^1 + d_{22} {\Omega}^2 - d_{12}{\Omega}^3 = 0  \]
that must be compared to the linearized Riemann operator for the euclidean metric leading to:   \\
\[   {\Omega}_{11} \equiv  2 {\xi}^1_1=0, {\Omega}_{12} \equiv {\xi}^1_2 + {\xi}^2_1=0, {\Omega}^{22} \equiv 2 {\xi}^2_2=0      \,\,\, \Rightarrow \,\,\, 
        d_{11}{\Omega}_{22} + d_{22}{\Omega}_{11} - 2 d_{12}{\Omega}_{12}=0  \]
for the euclidean metric.
In a more contructive way, we have ${\omega}^5 = {\partial}_1{\omega}^1 + {\omega}^1 \,\, {\omega}^4$ where ${\omega}^4$ is given in a rational way by the formula:  \\
\[   {\omega}^3 {\partial}_2{\omega}^2 - {\partial}_1{\omega}^3 + {\omega}^2 {\partial}_2{\omega}^3 - 
{\omega}^2 {\omega}^3 {\partial}_1{\omega}^1 +
 {\omega}^3 (1 - {\omega}^1 {\omega}^2) {\omega}^4=0             \]
Accordingly, if we do want to solve the equivalence problem ${\Phi}^1 = {\bar{\omega}}^1, {\Phi}^2 = {\bar{\omega}}^2, {\Phi}^3= {\bar{\omega}}^3$, {\it we must know the Vessiot structure equation}. As for the Vessiot structure constant $c$, it must be the same at first sight because the Vessiot structure equations are invariant under any diffeomorphism. However, we have to take into account the fact that two sections $\omega$ and ${\bar{\omega}}$ of ${\cal{F}}$ may give the same infinitesimal Lie equations. For example, in the present situation, we must have:   \\
\[       {\bar{\omega}}^1= {\omega}^1, \hspace{1cm} {\bar{\omega}}^2= {\omega}^2, \hspace{1cm}
{\bar{\omega}}^3 = a {\omega}^3 \hspace{1cm} \Rightarrow \hspace{1cm} \bar{c}= c/a   \]
where $a\neq 0$ is the parameter of the multiplicative group of the real line. It follows that we have $c=0 \Rightarrow \bar{c}=0$ both with $c\neq 0 \Rightarrow \bar{c}\neq 0$ and it just remains to exhibit such situations.  \\
In the present situation, we have ${\omega}^1=0, {\omega}^2=0, {\omega}^3=1 \Rightarrow c=0$  \\ 
However, the new pseudogroup:   \\
\[  \bar{\Gamma}=\{ y^1=\frac{ax^1 +b}{cx^1 + d}, y^2= \frac{ax^2 +b}{cx^2 + d} \mid a,b,c,d= cst  \}  \] 
is easily seen to be provided by the new specialization:  \\
\[ {\bar{\omega}}^1=0, {\bar{\omega}}^2=0, {\bar{\omega}}^3= 1/((x^2 - x^1)^2 \Rightarrow \bar{c}= - 2  \]  
leading to the new Lie form:  \\
\[    y^1_2=0, \,\,\, y^2_1=0,  \,\,\, \frac{1}{(y^2 - y^1)^2}y^1_1 y^2_2 = \frac{1}{(x^2 - x^1)^2} \]
It follows that the equivalence problem $y^1_2/y^1_1=0, y^2_1/y^2_2=0, y^1_1y^2_2= 1/ (x^2-x^1)^2$ cannot be solved. \\
Indeed, we should get $y^1=f(x^1), y^2=g(x^2)$ and thus $y^1=f(x^1), y^2=g(x^2)$ with ${\partial}_1f(x^1){\partial}_2g(x^2)= 1/(x^2 - x^1)^2$. Inverting the formula and setting $x^1=x^2=x$, we should conclude that ${\partial}_x f(x)=0$ or ${\partial}_x g(x)=0$ but this is impossible. \\
We could also notice that $0 = ac =\bar{c} =  - 2$ for a certain $a\neq 0$ and a contradiction.  \\
Needless to say that no classical tool can produce anyone of these results.   \\

\newpage

\noindent
{\bf 3) RIEMANN STRUCTURE}   \\

In such a way that this section could be even more striking, we shall copy almost " word by word " the procedure of the preceding section. \\
With $m=n=2$, let us consider the Lie group of isometries $y=Ax+B$ where $A$ is an othogonal matrix for the eucldean metric $\omega=(dx^1)^2 + (dx^2)^2$ as an algebraic Lie pseudogroup $\Gamma$. It is easy to exhibit the corresponding first order system ${\cal{R}}_1$ of finite Lie equations in Lie form by introducing the three generating differential invariants ${\Phi}_{ij}\equiv {\omega}_{kl}(y)y^k_iy^l_j $ and the corresponding Lie form:  \\
\[  {\Phi}_{11} \equiv (y^1_1)^2 + (y^2_1)^2=1, \,\,\,   {\Phi}_{22} \equiv (y^1_2)^2 + (y^2_2)^2=1, \,\,\,  {\Phi}_{12} \equiv  y^1_1y^1_2 + y^2_1y^2_2=0  \]
The details of the following tricky computations, first done in 1978 ([8]), have been improved in 2016 ([14]) and are again revisited in this paper. Its major interest is to work out the two Vessiot structure constant existing but within a tensorial framework now. By this way, we prove that the well known formal integrability result found by L. P. Eisenhart in 1926 ([2]) on the constant Riemannian curvature condition is only a very paticular case of the Vessiot structure equations found by E. Vessiot more than twenty years before ([21]). \\
First of all, we notice that this system is finite type with a vanishing second order symbol and is thus formally integrable but not involutive. For this, we may introduce the six generating differential invariants ${\Phi}^k_{ij}$ obtained after one prolongation and transforming like the six Christoffel symbols of Riemannian geometry:  \\
\[   {\gamma}^k_{ij}=\frac{1}{2} {\omega}^{kr}({\partial}_i{\omega}_{rj} + {\partial}_j{\omega}_{ir} - {\partial}_r{\omega}_{ij}) = {\gamma}^k_{ji} \,\,\,  \Rightarrow \,\,\,  2 {\Phi}_{kr}{\Phi}^r_{ij} \equiv d_i {\Phi}_{rj} + d_j{\Phi}_{ir} - d_r {\Phi}_{ij}    \]
by ntroducing the inverse matrix of $\omega$. We notice that:   \\
\[ det({\Phi}_{ij}) = {\Phi}_{11}{\Phi}_{22} - ({\Phi}_{12})^2 \equiv (y^1_1 y^2_2 - y^1_2 y^2_1)^2\neq 0  \]
is changing like the square of the Jacobian $\Delta$ and we have ${\gamma}^r_{ri}=\frac{1}{2} {\omega}^{rs}{\partial}_i{\omega}_{rs}={\partial}_i(det(\omega))^{\frac{1}{2}}$.\\
Looking to the way these invariants are transformed under an arbitrary change $ \bar{x}=\varphi (x)$ of local coordinates, we obtain $u_{ij}={\partial}_i{\varphi}^k{\partial}_j{\varphi}^l {\bar{u}}_{kl}$ for describing the natural fiber bundle ${\cal{F}}= S_2T^*$ with local coordinates $(x^1,x^2; u_{11},u_{22},u_{12})$ and section $\omega = ({\omega}_{11},{\omega}_{22}, {\omega}_{12})$ becoming $(1,1,0)$ with our choice and $det(\omega)\neq 0$. Passing to the infinitesimal point of view we obtain the first order system $R_1 \subset J_1(T)$ in the Medolaghi form, also called Killing system, with jet notation $\Omega \equiv L({\xi}_1) \omega=0$ as in ([8],[11])  in dimension $n=2$, using capital letters for the linearization:  \\
\[   {\Omega}_{ij} \equiv {\omega}_{rj}(x){\xi}^r_i + {\omega}_{ir}(x) {\xi}^r_j + {\xi}^r {\partial}_r {\omega}_{ij}(x) = 0   \]
 that we can extend to six {\it intermediate} new second order equations equations, namely:  \\
\[  {\Gamma}^k_{ij} \equiv {\xi}^k_{ij} + {\gamma}^k_{rj}(x) {\xi}^r_i + {\gamma}^k_{ir}(x) {\xi}^r_j  - 
{\gamma}^r_{ij}(x) {\xi}^k_r + {\xi}^r {\partial}_r {\gamma}^k_{ij}(x) =0  \, \Rightarrow \, 
{\Gamma}^r_{ri}\equiv {\xi}^r_{ri} + {\gamma}^s_{sr}{\xi}^r_i  + {\xi}^r {\partial}_r {\gamma }^s_{si} = 0\]
It is only now that we have to use specific concepts of Riemannian geometry, namely the {\it Riemann} and {\it Ricci} tensors, both with the (second) {\it Bianchi} identities:  \\
\[     {\rho}^k_{l,ij} \equiv  {\partial}_i {\gamma}^k_{lj} - {\partial}_j{\gamma}^k_{li} + {\gamma}^r_{lj} {\gamma}^k_{ri} - {\gamma}^r_{li} {\gamma}^k_{rj} \,\,\, \Rightarrow \,\,\,   {\rho}_{ij} = {\rho}^r_{i,rj} \,\,\,  \Rightarrow \,\,\, 
{\varphi}_{ij}={\rho}_{ij} - {\rho}_{ji}= {\rho}^r_{r,ij}= {\partial}_i{\gamma}^r_{rj} - {\partial}_j{\gamma}^r_{ri} \]
Accordingly, {\it if we use the Christoffel symbols} $\gamma$ {\it independently of} ${\omega}$, we may have ${\rho}_{ij} \neq  {\rho}_{ji}$ ({\it care}) and we have the following lemma:  \\

\noindent
{\bf LEMMA 3.1}: The $\gamma$ alone are the geometric objects for the affine group $y=Ax + B$ with $n(n+1)$ parameters, which are a section of an affine natural bundle modeled on $S_2T^* \otimes T$. Introducing the Spencer $\delta$-map, the $n^2(n+1)/2$  generating differential invariants $({\Phi}^k_{ij})$ allow to describe a vanishing Riemann tensor in the Janet sequence for the Lie operator $T \rightarrow S_2T^* \otimes T : \xi \rightarrow \Gamma \equiv {\cal{L}}(\xi)\gamma $:  
\[    0 \rightarrow \Theta \rightarrow T \rightarrow S_2T^*\otimes T \rightarrow F_1 \rightarrow ... \]
where $F_1$ is defined by the short exact sequence $0 \rightarrow F_1 \rightarrow {\wedge}^2 T^* \otimes T^* \otimes T \stackrel{\delta}{\longrightarrow} {\wedge}^3 T^* \otimes T \rightarrow 0 $.  \\ 

\noindent
{\it Proof}: As $g_2=0$ and thus $g_3=0$, we have the commutative and exact diagram:  \\
\[   \begin{array}{rcccccccl}
    &   &    &   0  & &  0  &  &  &     \\
     &  &    &   \downarrow  &  &  \downarrow  &   &  &  \\
   &  0  & \rightarrow &  S_3T^*\otimes T  &  \rightarrow &  T^* \otimes S_2T^*\otimes T^* &  \rightarrow & F_1 & \rightarrow 0   \\
   &   &   &  \hspace{3mm}  \downarrow \delta &   &  \parallel &  &   \\
 &  0  & \rightarrow & T^* \otimes  S_2T^*\otimes T  &  =  &  T^* \otimes S_2T^*\otimes T &  \rightarrow & 0 &   \\          
&   &   & \hspace{3mm} \downarrow \delta  &   &  \downarrow &  &   \\
 &   &   &  \underline{{\wedge}^2T^* \otimes T^* \otimes T} &  &   0   &  &    &   \\
   &   &   & \hspace{3mm}  \downarrow \delta  &   &     &  &   &  \\
   &  &  &  {\wedge}^3 T^* \otimes T &   &    &  &  &  \\
   &  &   &   \downarrow   &  &  &  &     \\
   &   &  &   0  &  &     &   &   &
   \end{array}    \]
   and deduce from a snake chase that $F_1 \simeq \delta (T^* \otimes S_2T^* \otimes T)$.  \\
Now, from the transformation rules of $\gamma$, we deduce that $y^k_r {\Phi}^r_{ij} = y^k_ {ij}$ because $\gamma=0$ for the general affine group. Differentiating formally and substituting, we obtain the Vessiot structure equations: 
\[      d_i{\Phi}^k_{lj} - d_j {\Phi}^k_{li}+ {\Phi}^r_{lj} {\Phi}^k_{ri} - {\Phi}^r_{li}{\Phi}^k_{rj} = 0  \]
In dimension $2$, we have $F_1 \simeq {\wedge}^2T^* \otimes T^* \otimes T$ with $dim(F_1)=4$ like in the previous section.  \\
\hspace*{13cm}   Q.E.D.   \\

The linearization provides:  \\
\[  R^k_{l,ij} \equiv   - {\rho}^r_{l,ij} {\xi}^k_r + {\rho}^k_{r,ij} {\xi}^r_l + {\rho}^k_{l,rj} {\xi}^r_i + {\rho}^k_{l,ir} {\xi}^r_j + {\xi}^r {\partial}_r {\rho}^k_{l,ij} =0 \,\,  \Rightarrow \,\,  R_{ij} \equiv {\rho}_{rj}{\xi}^r_i + {\rho}_{ir}{\xi}^r_j + {\xi}^r {\partial}_r {\rho}_{ij} =0  \]
\[   F_{ij} \equiv R^r_{r,ij} \equiv {\varphi}_{rj}{\xi}^r_i + {\varphi}_{ir} {\xi}^r_j + {\xi}^r {\partial}_r {\varphi}_{ij}=0  \]

In the specific dimension $n=2$ considered, we have by chance the simplified formulas: \\
\[  {\rho}_{11}={\rho}^r_{1,r1}={\rho}^2_{1,21}, \,\, {\rho}_{12}={\rho}^r_{1,r2}= {\rho}^1_{1,12}, \,\, {\rho}_{21}={\rho}^r_{2,r1} = {\rho}^2_{2,21}, \,\, 
{\rho}_{22}={\rho}^r_{2,r2}={\rho}^1_{2,12}  \]   
and thus the {\it specific isomorphism} ${\wedge}^2 T^* \otimes T^* \otimes T \simeq {\wedge}^2 T^* \oplus S_2T^*$ {\it only for } $n=2$, defined by: \\
\[     ({\rho}^k_{l,ij}) \rightarrow ( {\rho}_{ij}) \rightarrow (\frac{1}{2}({\rho}_{ij} - {\rho}_{ji}), \frac{1}{2}({\rho}_{ij} + {\rho}_{ji}))  \]
by counting the dimensions with $1 \times 2 \times 2 =4 =1 +3$ while using the canonical splitting of the short exact $\delta$-sequence: 
\[  0 \rightarrow S_2T^* \stackrel{\delta}{\longrightarrow} T^* \otimes T^* \stackrel{\delta}{\longrightarrow}{\wedge}^2T^* \rightarrow 0  \]
Now, we have proved in many books ([8],[11],[14],[15]) or papers that, for any dimension $n$, two sections $\omega$ and $\bar{\omega}$ provides the same system of infinitesimal Lie equations if and only if $\bar{\omega} = a \omega$ for $a\neq 0$ the parameter of the multiplicative group of the real line. It follows that we have necessarily a first Vessiot constant $c_1$ in such a way that:  \\
\[   \frac{1}{2} ((\rho_{ij} + {\rho}_{ji}) = c_1 {\omega}_{ij}  \]
However, after linearization, we also obtain:  \\
\[  2 det(\omega) ({\xi}^1_1 + {\xi}^2_2) + {\xi}^r   {\partial}_r det(\omega) = 0  \]
and obtain therefore a second Vessiot structure constant $c_2$ such that we have (Compare [8] to the computer algebra result found in [5]):   \\
\[   \frac{1}{2}({\rho}_{12} - {\rho}_{21}) = \frac{1}{2} {\varphi}_{12}= c_2 (det(\omega))^{ \frac{1}{2}}   \] 
Finally, we have to take into account that we do want second order {\it integrability conditions} for the metric $\omega$, that is we must eliminate $\gamma $ by using the Levi-Civita isomorphism $(\omega,\gamma) \simeq j_1(\omega)$. Then, it is well known that ${\varphi}_{ij}=0$ and we must thus have $c_2=0$. Also, using the same diagram as in the previous section, we must have only one second order integrability condition which is indeed the constant curvature conditon expressed by means of the Ricci tensor which is now symmetric.  \\

\noindent
{\bf REMARK 3.2}: Let us prove that there is almost no difference with the example presented in the preceding section (not quoted in [5]), even though the background group is quite different. For this, let us introduce the different generating differential invariants and new Lie form:  \\
\[    {\Phi}_{11}= 2{\Phi}^2 {\Phi}^3\equiv 2 y^1_1 y^2_1=0, \,\, {\Phi}_{22}= 2{\Phi}^1 {\Phi}^3 \equiv 2 y^1_2 y^2_2=0, \,\, 
{\Phi}_{12} = {\Phi}^3 +{\Phi}^1{\Phi}^2{\Phi}^3 \equiv y^1_1 y^2_2 + y^1_2 y^2_1 = 1  \]
in such a way that:
\[   {\omega}_{11} = 2 {\omega}_2 {\omega}_3, {\omega}_{22}= 2 {\omega}_1 {\omega}_3, {\omega}_{12}={\omega}_3 (1 + {\omega}_1 {\omega}_2) \,\, 
\Rightarrow \,\,  d_{11}{\Omega}_{22} + d_{22} {\Omega}_{11} - 2 d_{12} {\Omega}_{12} = 0  \]
We now notice that ${\Phi}_{11}{\Phi}_{22} - ({\Phi}_{12})^2 = - {\Delta}^2 $ and obtain again ${\Phi}_{ij}\equiv{\omega}_{kl}y^k_i y^l_j$ as before but now with the strange metric ${\omega}_{11}=0, {\omega}_{22}=0, {\omega}_{12}=1$ in such a way that $det(\omega)= - 1 < 0$ contrary to the previous example where ${\bar{\omega}}_{11}=1,{\bar{\omega}}_{22}=1, {\bar{\omega}}_{12}=0$ leading to $det(\bar{\omega})=1> 0$. It follows that, even if the two examples can be treated similarly, there is no way to pass from one example to the other by an invertible transformation. Indeed, if we coulf find $y=f(x)$ such that ${\omega}_{kl}(x){\partial}_if^k(x){\partial}_jf^l(x) = {\bar{\omega}}_{ij}(x)$, taking the determinants, we should obtain $  det(\omega) {\Delta}^2=det(\bar{\omega})=1$ and the contradiction ${\Delta}^2 = - 1$.  \\

\noindent
{\bf REMARK 3.3}: When $n=3$, the Lie pseudogroup of transformations preserving the contact $1$-form $\alpha = dx^1 - x^3 dx^2$ is unimodular because it also preserves $\beta =d \alpha= dx^2 \wedge dx^3$ and thus the volume form $\alpha \wedge \beta = dx^1 \wedge dx^2 \wedge dx^3$. Hence, we can consider the geometric object $\omega =(\alpha, \beta)$ with $\alpha \wedge \beta \neq 0$ which is a meaningful condition like $det(\omega)\neq 0$ previously. In this case, the Vessiot structure equations are 
$d\alpha = c' \beta, d \beta = c'' \alpha \wedge \beta $ with two Vessiot structure constants (See [14] for details). Closing this exterior system, we obtain:  \\
\[   0 =d (d\alpha )= c' d\beta =c'c'' \alpha \wedge \beta \hspace{1cm} \Rightarrow \hspace{1cm} c'c''= 0  \]
which is quite unusual. Hence, one of the two constants must vanish like before but for a completely different reason. More generally, we refer the reader to ([16]) for the study of how the extension modules in homological algebra may depend on the Vessiot structure contants.   \\

\noindent
{\bf 4) CONCLUSION}   \\

According to the italian mathematician U. Amaldi in 1907, at the beginning of the last century, only two frenchmen, namely  E. Cartan and E. Vessiot, were knowing and understanding the work of S. Lie on the infinite groups of transformations, now called Lie pseudogroups . However, the respective {\it Cartan structure equations} and {\it Vessiot structure equations} that they developed about at the same time were so different that Amaldi said that " only the future should say which of the two should be the most important one ". Unhappily and mainly for private reasons that we have explained, Cartan and followers never told that there could be another way superseding their approach. Also Cartan never told to A. Einstein in his letters of 1930 on absolute parallelism about the work of M. Janet in 1920 on systems of partial differential equations ([3]). Then, D.C. Spencer, largely ignoring these tentatives, created around 1970 new tools for the formal study of systems of partial differential equations ([20]), in particular the ones allowing to define Lie pseudogroups ([4]). As a matter of fact, the nonlinear Spencer sequences are superseding the Cartan structure equations because they are able to quotient the results down to the ground manifold while Cartan was doing exterior calculus on jet bundles but this result is largely unknown today. As a byproduct, the work of Vessiot is still almost totally unknown after more than a century. It is only in 2008 that an interesting PhD thesis has been done by a german student who wanted to use computer algebra for exhibiting the Vessiot structure constants in the Riemannian case ([5]). However, it is clear that the student was lucky to be able to use the corresponding tensorial framework. The purpose of the present paper has been to illustrate and justify these results by means of a general constructive procedure, proving that things are not so simple for Riemannian or contact structures. We have already proved that the case of a conformal structure is even worst as it highly depends on the dimension of the ground manifold ([13][14],[16],[18]).   \\

\newpage

\noindent
{\bf REFERENCES}  \\

\noindent
[1] Cartan, E.: Sur la Structure des Groupes Infinis de Transformations, Ann. Ec. Norm. Sup., 21 (1904) 153-206.  \\
\noindent
[2] Eisenhart, L.P.: Riemannian Geometry, Princeton University Press, Princeton (1926).\\
\noindent
[3] Janet, M.: Sur les Syst\`{e}mes aux D\'{e}riv\'{e}es Partielles. Journal de Math\'{e}matique, 8 (1920) 65-151. \\
\noindent 
[4] Kumpera, A., Spencer, D.C.: Lie Equations, Ann. Math. Studies 73, Princeton University Press, Princeton (1972).\\
\noindent
[5] Lorenz, A.: On Local Integrability Conditions of Jets Groupoids, Acta Applicandae Mathematicae, 101 (2008) 205-213. \\
https://doi.org/10.1007/s10440-008-9193-7   \\
https://arxiv.org/abs/0708.1419v1  \\
\noindent
[6] Medolaghi, P.: Sulla Theoria dei Gruppi Infiniti Continui, Ann. Math. Pura Appl., 25 (1897) 179-218.  \\
\noindent
[7] Northcott, D.G.: An Introduction to Homological Algebra, Cambridge university Press, Cambridge (1966).  \\
\noindent
[8] Pommaret, J.-F.: Systems of Partial Differential Equations and Lie Pseudogroups. Gordon and Breach, New York (1978). Russian translation: MIR, Moscow (1983).\\
\noindent
[9] Pommaret, J.-F.: Differential Galois Theory. Gordon and Breach, New York (1983). \\
\noindent
[10] Pommaret, J.-F.: Lie Pseudogroups and Mechanics, Gordon and Breach, New York (1988). \\
\noindent
[11] Pommaret, J.-F.: Partial Differential Equations and Group Theory, Kluwer, Dordrecht (1994). \\
https://doi.org/10.1007/978-94-017-2539-2    \\
\noindent
[12] Pommaret, J.-F.: The Mathematical Foundations of General Relativity Revisited. Journal of Modern Physics, 4 (2013) 223-239. \\
 https://doi.org/10.4236/jmp.2013.48A022   \\
\noindent
[13] Pommaret, J.-F.: Airy, Beltrami, Maxwell, Einstein and Lanczos Potentials revisited. Journal of Modern Physics, 7 (2016) 699-728. \\
https://doi.org/10.4236/jmp.2016.77068   \\
\noindent
[14] Pommaret, J.-F.: Deformation Theory of Algebraic and Geometric Structures, Lambert Academic Publisher (LAP), Saarbrucken, Germany (2016). A short summary can be found in "Topics in Invariant Theory ", S\'{e}minaire P. Dubreil/M.-P. Malliavin, Springer Lecture Notes in Mathematics, 1478, pp 244-254. (1990). https://arxiv.org/abs/1207.1964  \\
\noindent
[15] Pommaret, J.-F.: New Mathematical Methods for Physics. Mathematical Physics Books, Nova Science Publishers, New York (2018). \\
\noindent
[16] Pommaret, J.-F.: Differential Homological Algebra and General Relativity. Journal of Modern Physics, 10 (2019) 1454-1486. 
https://doi.org/10.4236/jmp.2019.1012097   \\
\noindent
[17] Pommaret, J.-F.: A Mathematical Comparison of the Schwarzschild and Kerr Metrics, Journal of Modern Physics, 11, 10 (2020) 1672-1710.  \\
https://doi.org/10.4236/jmp.2020.1110104  \\
\noindent
[18] Pommaret, J.-F. (2020) The Conformal Group Revisited, https://arxiv.org/abs/2006.03449 .  \\
\noindent
[19] Rotman, J.J.: An Introduction to Homological Algebra. Academic Press (1979), Springer, New york (2009).   \\
\noindent
[20] Spencer, D.C.: Overdetermined Systems of Partial Differential Equations. Bulletin of the American Mathematical Society, 75 (1965) 1-114.\\
\noindent
[21] Vessiot, E.: Sur la Th\'{e}orie des Groupes Infinis. Annales de l'Ecole Normale Sup\'{e}rieure, 20 (1903) 411-451.   \\ 
\noindent

\end{document}